\documentclass[11pt]{amsart}
\usepackage{amssymb, adjustbox, amsbsy}
\usepackage{array, longtable, booktabs, tabularx}
\usepackage{geometry}

\newcolumntype{Y}{>{\centering\arraybackslash}X}

\usepackage{amsfonts, amscd}
\numberwithin{equation}{section}

\usepackage{bm}
\usepackage{verbatim}
\usepackage{graphicx}
\graphicspath{{images/}}
\usepackage{tikz, tikz-cd}
\usepackage{subcaption}
\usepackage{listings}
\usepackage{subfiles}
\usepackage{mathtools}
\usepackage{comment}
\usepackage{enumitem}
\usepackage[all]{xy}
\usepackage{hyperref}
\hypersetup{
    colorlinks=true,
    citecolor=red,
    linkcolor=blue,
    filecolor=magenta,      
    urlcolor=red,
}

\usepackage[textsize=footnotesize,bordercolor=red,linecolor=red,backgroundcolor=red!15]{todonotes}

\newcommand{\CC}{\mathbb{C}}

\newcommand{\Ker}{\operatorname{Ker}}

\newtheorem{thm}{Theorem}[section]

\theoremstyle{definition}

\newtheorem{ques}[thm]{Question}
\newtheorem{rem}[thm]{Remark}

\newtheorem{ex}[thm]{Example}

\tikzstyle{wbullet}=[circle, draw=black, fill=white, thick, inner sep=2pt, minimum size=1.5mm]
\tikzstyle{bbullet}=[circle, draw=black, fill=black, inner sep=2pt, minimum size=1.5mm]

\begin{document}

\title{An example of a very non-movable effective divisor}

\author{Jihao Liu}

\address{Department of Mathematics, Peking University, No. 5 Yiheyuan Road, Haidian District, Beijing 100871, China}
\address{Beijing International Center for Mathematical Research, Peking University, No. 5 Yiheyuan Road, Haidian District, Beijing 100871, China}
\email{liujihao@math.pku.edu.cn}
\date{\today}

\begin{abstract}
We give a negative answer to a question of Ciliberto, Knutsen, Lesieutre, Lozovanu, Miranda, Mustopa, and Testa on effective divisors of positive self-intersection on smooth projective surfaces. The main result of this paper is obtained by generative AI, particularly Chatgpt 5.5 pro and the Rethlas system.
\end{abstract}

\maketitle

\section{Introduction}\label{sec:introduction}

We work over the field of complex numbers $\mathbb C$.

Ciliberto, Knutsen, Lesieutre, Lozovanu, Miranda, Mustopa, and Testa \cite[Question~1]{Cil+17} asked the following question.

\begin{ques}\label{ques:knutsen-question}
Let $X$ be a smooth projective surface. Does there exist a constant $m_1 = m_1(X)$ such that if $D$ is any divisor with
$$h^0(X,\mathcal{O}_X(D))=1\quad \text{and}\quad D^2>0,$$ 
then one has $h^0(X,\mathcal{O}_X(m_1D))\geq 2$?
\end{ques}

As mentioned in \cite{Cil+17}, the question can be traced back to A. L. Knutsen during the Warsaw workshop \emph{Okounkov bodies and Nagata type Conjectures} held at the Banach Centre in September 2013 and also during the workshop \emph{Recent advances in linear series and Newton--Okounkov bodies} held in Padua in February 2015. We show that Question~\ref{ques:knutsen-question} has a negative answer by considering the following example.

\begin{ex}\label{ex:non-movable-divisors}
Let $E$ be a smooth elliptic curve. Let
$$Y:=E\times E,\quad F=\{0\}\times E,\quad \text{and}\quad G=E\times \{0\}.$$
For an odd integer $n\geq 3$, define
$$\varphi_n:Y\to E,\qquad (x,y)\mapsto nx+2y,$$
and set
$$\Gamma_n=\Ker(\varphi_n).$$
Since $\gcd(n,2)=1$, we may choose integers $a,b$ with $nb-2a=\pm 1$. Then the matrix
$$
\begin{pmatrix}
n&2\\
a&b
\end{pmatrix}
$$
defines an automorphism of $E\times E$, and hence $\Gamma_n$ is a connected smooth elliptic curve.

We first record the relevant intersections. We have
$$F\cap \Gamma_n=\{0\}\times E[2],\qquad G\cap \Gamma_n=E[n]\times \{0\}.$$
These intersections are scheme-theoretic, since multiplication by $2$ and by $n$ are \'etale in characteristic zero. Thus
$$F\cdot \Gamma_n=4,\qquad G\cdot \Gamma_n=n^2.$$
Moreover $F^2=G^2=0$, and $\Gamma_n^2=0$ by adjunction, since $\Gamma_n$ is a smooth elliptic curve on an abelian surface. Therefore, if
$$A_n=F+\Gamma_n,$$
then
$$A_n^2=F^2+2F\cdot \Gamma_n+\Gamma_n^2=8.$$

The set $F\cap \Gamma_n=\{0\}\times E[2]$ is independent of $n$. At any point of this set, the tangent spaces are
$$T(F)=\{(0,v)\},\qquad T(\Gamma_n)=\{(u,v)\mid nu+2v=0\}.$$
Their intersection is zero. Hence $F$ and $\Gamma_n$ meet transversely at these four points, and $A_n$ has ordinary nodes there. Fix three of these points once and for all:
$$S=\{p_1,p_2,p_3\}\subset \{0\}\times E[2].$$

Now put
$$R=F+G.$$
The line bundle
$$\mathcal O_Y(2R)\simeq \operatorname{pr}_1^*\mathcal O_E(2[0])\otimes \operatorname{pr}_2^*\mathcal O_E(2[0])$$
is base-point-free. Hence we can choose a smooth divisor
$$B\in |2R|$$
such that $B\cap S=\emptyset$. Let
$$f:X=\operatorname{Spec}_Y\bigl(\mathcal O_Y\oplus \mathcal O_Y(-R)\bigr)\to Y$$
be the double cover branched along $B$. Since $B$ is smooth, $X$ is smooth, and
$$f_*\mathcal O_X=\mathcal O_Y\oplus \mathcal O_Y(-R).$$
Also $R$ is ample, so $H^0(Y,\mathcal O_Y(-R))=0$. Thus
$$H^0(X,\mathcal O_X)=H^0(Y,\mathcal O_Y)=\CC,$$
and $X$ is connected, hence irreducible.

Since $B\cap S=\emptyset$, the morphism $f$ is \'etale over each $p_i$. Choose one point
$$q_i\in f^{-1}(p_i)$$
for each $i=1,2,3$, and let
$$\pi:\widetilde X\to X$$
be the blow-up of $q_1,q_2,q_3$, with exceptional curves $E_1,E_2,E_3$. Define
$$D_n=\pi^*f^*A_n-2(E_1+E_2+E_3).$$
Since $f$ is \'etale over the points $p_i$, the divisor $f^*A_n$ has an ordinary node at each $q_i$. Hence $D_n$ is effective: it is obtained by taking the strict transform at the three chosen nodes. Its self-intersection is
$$D_n^2=(f^*A_n)^2+4(E_1+E_2+E_3)^2=2A_n^2-12=4>0.$$

We claim that
$$h^0(\widetilde X,\mathcal O_{\widetilde X}(mD_n))=1$$
for every integer
$$1\leq m< n^2/4+1.$$
Let $Q=\{q_1,q_2,q_3\}$ and set $M=mA_n$. By the definition of $D_n$,
$$
H^0(\widetilde X,\mathcal O_{\widetilde X}(mD_n))
=H^0\bigl(X,\mathcal O_X(f^*M)\otimes I_Q^{2m}\bigr).
$$
The projection formula gives
$$
H^0(X,\mathcal O_X(f^*M))
=H^0(Y,\mathcal O_Y(M))\oplus H^0(Y,\mathcal O_Y(M-R)).
$$
We have
$$
(mA_n-R)\cdot \Gamma_n
=m(F+\Gamma_n)\cdot \Gamma_n-(F+G)\cdot \Gamma_n
=4(m-1)-n^2.
$$
If $1\leq m<n^2/4+1$, then this number is negative. Thus $mA_n-R$ is not nef. Since $Y$ is an abelian surface, $mA_n-R$ is not effective. Thus we have
$$
H^0(\widetilde X,\mathcal O_{\widetilde X}(mD_n))
\simeq H^0(Y,\mathcal O_Y(mA_n)\otimes I_S^{2m})
$$
when $1\leq m<n^2/4+1$.

It remains to compute the last space. Let
$$\rho:Y'\to Y$$
be the blow-up of $S$, with exceptional curves $e_1,e_2,e_3$, and define
$$L_n=\rho^*A_n-2(e_1+e_2+e_3).$$
Since the points of $S$ are transverse intersection points of $F$ and $\Gamma_n$, if $F'$ and $\Gamma_n'$ denote the strict transforms of $F$ and $\Gamma_n$ on $Y'$, then
$$L_n=F'+\Gamma_n'.$$
Moreover,
$$
(F')^2=-3,\qquad (\Gamma_n')^2=-3,\qquad\text{and}\qquad  F'\cdot \Gamma_n'=1.
$$
Thus
$$L_n\cdot F'=-2\qquad \text{and}\qquad L_n\cdot \Gamma_n'=-2.$$

We now show that $|mL_n|$ has exactly one member for every $m\geq 1$. Of course $m(F'+\Gamma_n')\in |mL_n|$, so $|mL_n|$ has at least one member. Pick any $C\in |mL_n|$. Since
$$C\cdot F'=mL_n\cdot F'=-2m<0,$$
$F'$ is an irreducible component of $C$. Write $C=F'+C_1$. Then
$$
C_1\cdot \Gamma_n'
=(mL_n-F')\cdot \Gamma_n'
=-2m-1<0,
$$
so $\Gamma_n'$ is a component of $C_1$. Hence
$$C=F'+\Gamma_n'+C_2,\qquad C_2\in |(m-1)L_n|.$$
Induction on $m$ gives
$$C=mF'+m\Gamma_n'.$$
Consequently,
$$h^0(Y',\mathcal O_{Y'}(mL_n))=1$$
for every $m\geq 1$. Since
$$mL_n=\rho^*(mA_n)-2m(e_1+e_2+e_3),$$
we also have
$$
H^0(Y',\mathcal O_{Y'}(mL_n))
\simeq H^0(Y,\mathcal O_Y(mA_n)\otimes I_S^{2m}).
$$
Combining the preceding equalities, we obtain
$$h^0\left(\widetilde X,\mathcal O_{\widetilde X}(mD_n)\right)=1$$
for all
$$1\leq m<n^2/4+1.$$
In particular $h^0\left(\widetilde X,\mathcal O_{\widetilde X}(D_n)\right)=1$ and $D_n^2=4>0$. Taking $n$ to be a sufficiently large odd integer, we obtain a negative answer to Question~\ref{ques:knutsen-question}. 
\end{ex}

\begin{rem}
We note that the divisors $D_n$ in our example are reducible. It is interesting to ask whether one can make $D$ prime or irreducible.
\end{rem}

\begin{rem}
The main result of this paper is obtained by generative AI by simply asking 2 questions to Chatgpt 5.5 pro without any mathematical inputs. Chatgpt 5.5 pro provided the right example. The author later put the example and its proof Chatgpt 5.5 pro created into the Rethlas system's verifier. The verifier alerted that Chatgpt 5.5 pro's proof on the part that $M-R$ has no sections is wrong. The author then ran the Rethlas system based on Chatgpt 5.5 pro's output and obtained the correct correction, which ends up with the current version. See \cite{Ju+26} for a detailed introduction to the Rethlas system.

Due to the limitation of generative AI, it is possible that we have missed some related references in the literature, and we welcome any comments from experts.
\end{rem}

\subsection*{Acknowledgements}
The author was partially supported by the National Key R\&D Program of China \#\allowbreak 2024YFA1014400. The author would like to thank the Rethlas team, namely Haocheng Ju, Jiedong Jiang, Shurui Liu, Guoxiong Gao, Yuefeng Wang, Zeming Sun, Bin Wu, Liang Xiao, and Bin Dong, for their contributions to the development of Rethlas and its customized version used for the problem studied in this paper. The author would like to thank Kaiyuan Gu, Ruicheng Hu, and Sheng Qin for assistance with the verification of an earlier blueprint of this paper. The author would like to thank Ruochuan Liu and Gang Tian for constant support and encouragement.

\end{document}